\documentclass[11pt]{amsart} \usepackage{amssymb,latexsym,amsmath,amsthm}
\usepackage{graphicx}

\theoremstyle{plain} \newcounter{Theorem} \newtheorem{theorem}[Theorem]{Theorem}
\newtheorem{lemma}[Theorem]{Lemma} \newtheorem{proposition}[Theorem]{Proposition}
\newtheorem{corollary}[Theorem]{Corollary} \theoremstyle{remark} 
\newcounter{Remark} \newcounter{Example} 
\newtheorem{remark}[Remark]{Remark} \newtheorem{example}[Example]{Example}
\newcommand{\Rone}{\mathbb R\,}  
\newcommand{\None}{\mathbb N\,}  
\newcommand{\hs}{\mathcal H\,} \newcommand{\bh}{\mathcal{B}(\hs)} \newcommand{\nn}{\nonumber}
\newcommand{\bmtx}{\left[ \begin{matrix}}
\newcommand{\emtx}{\end{matrix} \right]}
\renewcommand{\L}{\lambda}
\newcommand{\A}{\alpha}
\newcommand{\D}{\delta}
\newcommand{\G}{\gamma}
\newcommand{\E}{\epsilon}
\newcommand{\txt}[1]{\;\mbox{#1}\;}
\newcommand{\tx}[1]{\,\mbox{#1}\,}
\newcommand{\noi}{\noindent}

\begin{document}

\title[Sums of Rank-One Operators]{Rank-One Decomposition of Operators and Construction of Frames}
 \author{Keri A.  Kornelson} \address{Texas A\&M University} \email{keri@math.tamu.edu}
 \author{David R.  Larson} \address{Texas A\&M University}
 \email{larson@math.tamu.edu}
 \keywords{frames,tight frames, rank-one operators, positive operators, rank-one decomposition}
\subjclass[2000]{Primary 42C15, 47N40; Secondary 47C05, 46B28}

\thanks{The research of the first author was supported in part by an NSF VIGRE postdoctoral fellowship.  The research of the second author was supported in part by a grant from the NSF}

\begin{abstract} The construction of frames for a Hilbert space $\hs$ can be equated to the
decomposition of the frame operator as a sum of positive operators having rank one.  This
realization provides a different approach to questions regarding frames with particular properties
and motivates our results.   We find a necessary and
sufficient condition under which any positive finite-rank operator $B$ can be expressed as a sum of
rank-one operators with norms specified by a sequence of positive numbers $\{c_i\}$. Equivalently,
this result proves the existence of a frame with $B$ as it's frame operator and with vector norms
$\{\sqrt{c_i}\}$.  We further prove that, given a non-compact positive operator $B$ on an
infinite dimensional separable real or complex Hilbert space, and given an infinite sequence
$\{c_i\}$ of positive real numbers which has infinite sum and which has supremum strictly less than
the essential norm of $B$, there is a sequence of rank-one positive operators, with norms given by
$\{c_i\}$, which sum to $B$ in the strong operator topology.

These results generalize results by Casazza, Kova\v{c}evi\'{c}, Leon, and
Tremain, in which the operator is a scalar multiple of the identity operator (or equivalently
the frame is a tight frame), and also results by Dykema, Freeman, Kornelson, Larson, Ordower, and
Weber in which $\{c_i\}$ is a constant sequence.

\end{abstract} \maketitle

\section*{Introduction}

The existence and characterization of frames with a variety of additional properties is an area
of active research which has resulted in recent papers including \cite{BF}, \cite{CKLT}, and
\cite{DFK}.  These questions can also be phrased in terms of the expression of
positive operators as sums of rank-one positive operators.  Therefore, although motivated by
applications in signal processing, the results have independent interest to operator theory.

Throughout this paper, $\hs$ will be a real or complex separable Hilbert space.  For $\mathbb{J}$
a real or countably infinite index set, the collection of vectors $\{x_j\}_{j \in \mathbb{J}}
\subset \hs$ is a \textit{frame} if there exist positive constants $D \geq C >0$ such that for every
$x \in \hs$: $$ C\|x\|^2 \leq \sum_{j \in \mathbb{J}}|\langle x, x_j \rangle |^2 \leq D\|x\|^2$$
A frame is called a \textit{tight frame} if $C=D$, and a \textit{Parseval frame} if $C=D=1$.

We use elementary tensor notation for a rank-one operator on $\hs$.  Given $u,v,x \in \hs$, the operator $u
\otimes v$ is defined by $ (u \otimes v)x = \langle x,v \rangle u \quad \txt{for} x \in \hs$. The
operator $u \otimes u$ is a projection if and only if $\|u\|=1$.

Given a frame  $\{x_i\}$ for $\hs$, the frame operator is defined to be the map $S: \hs \rightarrow \hs$ taking $x
\rightarrow \sum_i \langle x, x_i \rangle x_i$.  We can write $S$ in the form:

\begin{equation}\label{op}
 S = \sum_i x_i \otimes x_i
 \end{equation}
\noi where the convergence is in the strong operator topology (SOT). It is a well-known result that $\{x_i\}$ is a
tight frame with frame bound $\L$ if and only if $S = \L I$, where $I$ is the identity operator on $\hs$. (see
\cite{HL}) Therefore, finding a tight frame is equivalent to decomposing a scalar multiple of the identity as a
sum of rank-one operators.

Given $\hs$ a finite dimensional Hilbert space, and $\{c_i\}$ a sequence of positive real numbers, Casazza,
Kova\v{c}evi\'{c}, Leon, and Tremain found in \cite{CKLT} a necessary and sufficient condition (called the
Fundamental Frame Inequality) for the existence of a tight frame for $\hs$ with the norms of the vectors given by
$\{c_i\}$.  As stated above, this is equivalent to decomposing a scalar multiple of the identity as a sum of
rank-one positive operators with prescribed norms.  In Section \ref{one}, we generalize this result to the case
where the scalar is replaced with an arbitrary  positive invertible operator on $\hs$ (in fact, we state the
result for positive operators which are not necessarily invertible), thus obtaining the condition for the
existence of non-tight frames with specified norms.

In Section \ref{two}, $\hs$ is an infinite dimensional separable real or complex Hilbert space.  Recall that the
essential norm of an operator on $\hs$ is defined to be $\|B\|_{ess} = \inf\{\|B-K\|\}$ where $K$ is a compact
operator on $\hs$.  We prove that, given $B$ a non-compact positive operator on $\hs$ and given $\{c_i\}$ a
sequence of positive real numbers which sum to infinity and for which $\sup_i \{c_i\} < \|B\|_{ess}$, $B$ can be
expressed as a sum of rank-one operators having norms given by $\{c_i\}$, with convergence in the SOT. This
implies that every positive invertible operator with essential norm strictly greater than $\sup_i \{c_i\}$ is the
frame operator for a frame with prescribed norms $\{\sqrt{c_i}\}$. This result is the generalization of the
infinite-dimensional result in \cite{CKLT}, which restricts to tight frames. It also generalizes the results in
\cite{DFK} regarding ellipsoidal tight frames, in which the sequence $\{c_i\}$ is a constant sequence.

The authors would like to thank Eric Ricard for his helpful comments during the Linear Analysis
Workshop at Texas A\&M University in August, 2003.

\section{Finite Dimensions}\label{one}
The following lemma was shown to us by Marc Ordower.  It was a key step in the proof of Proposition
6 in \cite{DFK}.

\begin{lemma}\label{marc}  Let $B$ be a positive operator with rank $n$ and nonzero
eigenvalues $b_1 \geq b_2 \geq \cdots \geq b_n > 0$ and with a corresponding orthonormal set
of eigenvectors $\{e_i\}_{i=1}^n$.  If $c$ is a positive number with $b_j \geq c \geq
b_{j+1}$ for some $1 \leq j \leq n-1$, then there is a unit vector $x$ in span$\{e_j,
e_{j+1}\}$ such that $B-c(x \otimes x)$ is positive and has rank $n-1$.
\end{lemma}

\begin{proof}
If $A$ is a finite-rank self-adjoint operator and $c \geq \|A\|$ is a constant, then the
least eigenvalue of $A$ is given by $\|(A+cI)^{-1}\|^{-1}-c$.  It follows that the
function giving the least eigenvalue is continuous on the set of finite-rank self-adjoint
operators.  Let $A_x = B-c(x \otimes x)$ for $x$ an arbitrary unit vector in
span$\{e_j,e_{j+1}\}$.  The least eigenvalue of $A_x$ is nonnegative for $x=e_j$ and
nonpositive for $x=e_{j+1}$, hence span$\{e_j, e_{j+1}\}$ contains a unit vector $x$ for which the
least eigenvalue of $A_x$ is zero.  Therefore, $A_x$ is positive of rank $n-1$.

\end{proof}

\begin{remark}  Since the trace of $A_x$ is $\left(\sum_1^n b_i \right) - c\,$, and eigenvalues
other than $b_j$ and $b_{j+1}$ are unchanged, the remaining nonzero eigenvalue must be
$$\tilde{b}=b_j + b_{j+1}-c$$  Moreover, the new eigenvalue maintains the $j^{\tx{th}}$ position in
the ordering of all the eigenvalues: $$b_1 \geq \cdots \geq b_{j-1} \geq \tilde{b} \geq b_{j+2} \geq
\cdots \geq  b_n$$

\end{remark}

%%%%%%%%%%%%%%%%%%%%%%%%%%%%%%%%%%%%%%%%%%%%%%%%%
%%%%%%     Theorem                        %%%%%%%
%%%%%%%%%%%%%%%%%%%%%%%%%%%%%%%%%%%%%%%%%%%%%%%%%
\begin{theorem}\label{sums}  Let $B$ be a positive operator with rank $n$
and nonzero eigenvalues $b_1 \geq b_2 \geq \cdots \geq b_n >0$.  Let $\{c_i\}_{i=1}^k\,;\; k \geq n$ be a
nonincreasing sequence of positive numbers such that $\sum_{i=1}^k c_i = \sum_{j=1}^n b_j$.  There exist unit
vectors $\{x_i\}_{i=1}^k$ such that $$B = \sum_{i=1}^k c_i(x_i \otimes x_i)
$$
\noi if and only if for each $p$, $1 \leq p \leq n-1$,
\begin{equation}\label{ps} \sum_{i=1}^p c_i \leq \sum_{i=1}^p b_i
\end{equation}

\end{theorem}
%%%%%%%%%%%%%%%%%%%%%%%%%%%%%%%%%%%%%%%%%%%%%%%%%%

\noi \textit{Note: } When a positive invertible operator $B$ can be written in the form of this
theorem, we say that $B$ admits a \textit{rank-one decomposition corresponding to $\{c_i\}$}.

\begin{remark} The hypotheses in Theorem \ref{sums} reduces to exactly the Fundamental Frame
Inequality from \cite{CKLT} when $b_1=b_2= \cdots = b_n$. Thus, an immediate corollary of Theorem \ref{sums} is
the result from \cite{CKLT} that the Fundamental Frame Inequality is a necessary and sufficient condition for the
existence of a tight frame for $\Rone^n$ consisting of vectors with norms $\sqrt{c_i}\,;\,i=1, \ldots, k$. We
include this as Corollary \ref{pete}.

\end{remark}

\begin{proof}[Proof of Theorem \ref{sums}]

Assuming the property (\ref{ps}) holds, we use strong induction on the rank $n$ of the operator.  The case $n=1$ when
$B$ is already a rank-one operator is clear. Let $\{e_i\}_1^n$ be an orthonormal set of eigenvectors
for $b_1, b_2, \ldots, b_n$, respectively.

\noi \textit{Case 1:} If $b_1 \geq c_1 \geq b_n$, there exists an index $l\,;\, 1 \leq l \leq n-1$ such that $b_l
\geq c_1 \geq b_{l+1}$.  By Lemma \ref{marc}, there exists a unit vector $x$ in the span of $e_l$ and $e_{l+1}$
such that $B-c_1(x \otimes x)$ has rank $n-1$ and eigenvalues in nonincreasing order $b_1, \ldots, b_{l-1},
\tilde{b}, b_{l+2}, \ldots, b_n$, where $\tilde{b}=b_l+b_{l+1}-c_1$.  Let $x_1= x$.  It remains to check that the
inequalities (\ref{ps}) hold for the operator $B-c_1(x_1 \otimes x_1)$ and the sequence $\{c_2, \cdots, c_k\}$.
Clearly, for $1 \leq r \leq l$, we have $$ c_2 + \cdots +c_r \leq c_1 + \cdots + c_{r-1} \leq b_1 + \cdots
b_{r-1}$$ For $ l+1 \leq r \leq n-1$, \begin{eqnarray*} c_1 + \cdots + c_l +c_{l+1} +\cdots + c_r & \leq & b_1 +
\cdots + b_l + b_{l+1} + \cdots + b_r \\ \Rightarrow \; c_2 + \cdots + c_l +c_{l+1} +\cdots + c_r &\leq & b_1 +
\cdots + \left( b_l + b_{l+1}-c_1 \right) + \cdots + b_r \\ \Rightarrow \; c_2 + \cdots + c_l +c_{l+1} +\cdots +
c_r &\leq & b_1 + \cdots + \tilde{b} + \cdots + b_r
\end{eqnarray*}

Since we have reduced to an operator of rank $n-1$ which satisfies the hypotheses of the theorem, induction gives
the remaining elements of the rank-one decomposition of $B$ corresponding to the sequence $\{c_i\}_1^k$.

\noi \textit{Case 2:}  If $c_1 < b_n$, then we are unable to use Lemma \ref{marc}.  We select rank-one operators
to subtract from $B$ which preserve the rank and decrease the smallest eigenvalue until the Case 1
property is attained.

Let $p$ be the largest integer such that $c_1 + c_2 + \cdots + c_p < b_n$.  Note that Case 2 can only occur for
$k>n$.  Since the sums of each sequence are equal, $c_{p+1} + \cdots +c_k > b_1 + \cdots +b_{n-1}$,
but because $c_{p+1} \leq c_1 \leq b_n \leq b_{n-1}$, the sum on the left must have more terms.
Therefore, $k-p > n-1$, or alternatively, $1 \leq p \leq k-n$. Let $x_i = e_n$ for each $i=1, 2,
\ldots, p$. The operator $B-\sum_{i=1}^p c_i(x_i \otimes x_i)$ still has rank $n$ and the
eigenvalues in decreasing order are $b_1, b_2, \ldots, b_{n-1}, \tilde{b}$, where $\tilde{b} = b_n -
\sum_{i=1}^p c_i$.  By the selection of $p$, $c_{p+1}$ exceeds or equals $\tilde{b}$, the smallest
eigenvalue of $B-\sum_{i=1}^p c_i(x_i \otimes x_i)$. The method of Case 1 can now be applied to the
operator $B-\sum_{i=1}^p c_i(x_i \otimes x_i)$ and the remaining sequence $c_{p+1}, \ldots, c_k$.

Conversely, assume $B$ has a rank-one decomposition corresponding to $\{c_i\}$.  Given a fixed $j$ with
$1 \leq j \leq k$, define $P$ to be the orthogonal projection onto the span of $\{x_i\}_{i=1}^j$.
Clearly, rank $P \leq j$. For each $i$, define $P_i = x_i \otimes x_i$.  We then have: $$ PBP =
\sum_{i=1}^k c_i PP_iP \geq \sum_{i=1}^j c_i PP_iP = \sum_{i=1}^j c_iP_i$$ \noi and therefore
trace$(PBP) \geq \tx{trace}\left(\sum_{i=1}^j c_iP_i \right) = c_1 + c_2 + \cdots c_j$.  We next
show that trace$(PBP) \leq b_1 + b_2 + \cdots +b_n$, which will complete the proof.

Let $\{e_i\}_{i=1}^n$ be as above, and for notational purposes, set $b_{n+1} = 0$.  For $1 \leq i
\leq n$, define $Q_i = e_1 \otimes e_1 + \cdots + e_i \otimes e_i$.  A calculation verifies that $B$
can be written
$$B = \sum_{i=1}^n (b_i - b_{i+1})Q_i$$

\noi Then, $PBP = \sum_{i=1}^n (b_i - b_{i+1})PQ_iP$. Observe that trace$(PQ_iP) = \tx{trace} (PQ_i) \leq \min\{ \tx{rank} P, \tx{rank} Q_i\} \leq \min\{j,i\}$.  Therefore,
\begin{eqnarray*}
\tx{trace} (PBP) &=& \sum_{i=1}^n (b_i - b_{i+1}) \tx{trace} (PQ_iP)\\
&\leq& \sum_{i=1}^j (b_i - b_{i+1})i + j\sum_{i=j+1}^n (b_i - b_{i+1}) \\
&=& \left(b_1 + b_2 + \cdots + b_j -jb_{j+1}\right) + j\left(b_{j+1} - b_{n+1} \right)\\
 &=& b_1 + b_2 + \cdots + b_j
 \end{eqnarray*}

\end{proof}

Theorem \ref{sums} will be used in the proof of Proposition \ref{ident}.  To demonstrate the theorem, we have
included the following example for the case where $\hs = \Rone^2$.
%%%%%%%%%%%%%%%%%%%%%%%%%
%  Examples             %
%%%%%%%%%%%%%%%%%%%%%%%%%
\begin{example} Let $B = \bmtx 5&0\\0&4 \emtx$.  We find the rank-one decomposition of $B$
corresponding to $\{c_i\}_{i=1}^4 = 3,3,2,1$.

Since $c_1 < b_2$ (Case 2), we take $x_1 = \left[ \begin{matrix}1\\0 \end{matrix}\right]$.  This leaves $B-3(x_1
\otimes x_1) = \bmtx 5&0\\0&1 \emtx$ to decompose.  We now have $c_2$ between the two eigenvalues 5 and 1 (Case
1). Let $x(t) = \sqrt{1-t}\bmtx 1\\0 \emtx + \sqrt{t} \bmtx 0\\1 \emtx = \bmtx \sqrt{1-t} \\ \sqrt{t} \emtx$.  We
wish to find $t$ such that $\bmtx 5&0\\0&1 \emtx - 3(x(t) \otimes x(t))$ has rank one, i.e. the determinant is
zero. Straightforward calculation shows the solution $t = \frac{1}{6}$.  Take $x_2 = x(\frac{1}{6}) = \bmtx
\sqrt{\frac{5}{6}}\\ \sqrt{\frac{1}{6}} \emtx$.  The remainder  $\bmtx 5&0\\0&1 \emtx - 3(x_2 \otimes x_2) = \bmtx
\frac{5}{2}& -\frac{\sqrt{5}}{2}\\-\frac{\sqrt{5}}{2}& \frac{1}{2} \emtx $ has rank one and it's range is spanned
by the unit vector $z =\bmtx \sqrt{\frac{5}{6}}\\ -\sqrt{\frac{1}{6}} \emtx$.  We then must have $x_3 = x_4 = z$,
which completes the decomposition.
$$ 3\bmtx 0&0\\0&1\emtx + 3 \bmtx \frac{5}{6} & \frac{\sqrt{5}}{6} \\ \frac{\sqrt{5}}{6} & \frac{1}{6} \emtx +
(2+1)\bmtx \frac{5}{6} & -\frac{\sqrt{5}}{6} \\ -\frac{\sqrt{5}}{6} & \frac{1}{6} \emtx = \bmtx 5&0\\0&4 \emtx =
B$$
\end{example}

The next example demonstrates that a rank-one decomposition is impossible when property (\ref{ps}) does not hold.
 \begin{example}
The  operator $B= \bmtx 5&0&0\\ 0&2&0\\0&0&2 \emtx$ does not have a rank-one decomposition
corresponding to the sequence $\{c_i\}_{i=1}^3 = \{4,4,1\}$.  To see this, assume that we can find
unit vectors $x_1, x_2, x_3$ such that $B = 4(x_1 \otimes x_1) + 4 (x_2 \otimes x_2) + (x_3 \otimes
x_3)$.  Let $P$ be the projection onto the span of $\{x_1, x_2 \}$.  Then trace$(PBP) =
\sum_{i=1}^3 c_i \tx{trace} P(x_i \otimes x_i)P \geq 4+4 = 8$, because $P(x_i \otimes x_i)P = x_i \otimes x_i$, which has trace 1.

Given $P$ any projection with rank two, using the argument from the last paragraph of Theorem \ref{sums}, we have
 $\tx{trace}(PBP)$ less than or equal to 7, i.e. the sum of the largest two
eigenvalues.  The contradiction implies $B$ does not have a rank-one decomposition corresponding to $\{4, 4, 1\}$.

\end{example}

The following proposition lends insight into the underlying
geometry of tight frames.   We include it for independent interest.

%%%%%%%%%%%  Proposition   %%%%%%%%%%%%%%%%%%%
 \begin{proposition}[Alternate Version of $\Rone^2$ Case]\label{2d}  Let $B$ be a positive
invertible operator on $\Rone^2$, and express $B$ as  diag$(b_1, b_2)$, where $b_1
\geq b_2 >0$.  Let $\{c_i\}_{i=1}^k$ be positive constants such that $c_1 \geq c_2 \geq
\cdots \geq c_k$ and $$\sum_{i=1}^k c_i = b_1 + b_2 $$
If $c_1 \leq b_1$, then there
exist unit vectors $\{x_i\}_{i=1}^k$ such that:  $$ B = \sum_{i=1}^k c_i \left( x_i
\otimes x_i \right) $$

\end{proposition}
%%%%%%%%%%%%%%%%%%%%%%%%%%%%%%%%%%%%%%%%%%%%%%

 \begin{proof} Let $$x_i = \left[ \begin{matrix} \cos {\theta}_i \\ \sin {\theta}_i \end{matrix}
\right] \quad \text{so} \quad x_i \otimes x_i = \left[ \begin{matrix} \cos^2 {\theta}_i &
\cos {\theta}_i \sin {\theta}_i \\ \cos {\theta}_i \sin {\theta}_i & \sin^2 {\theta}_i
\end{matrix} \right] $$

Then let $B = \text{diag}(b_1, b_2)$ and $\{c_i\}_{i=1}^k$ a sequence of positive numbers
such that $\sum_1^kc_i = b_1 + b_2$ and $c_1 \geq c_2 \geq \cdots \geq c_k$. Then $B=
\sum_{i=1}^k c_i \left( x_i \otimes x_i \right) $ if and only if the following equations
hold:
\begin{eqnarray} c_1 \cos^2\theta_1 + c_2\cos^2\theta_2 + \cdots + c_k \cos^2\theta_k &=&
b_1 \nn \\ c_1 \sin^2\theta_1 + c_2\sin^2\theta_2 + \cdots + c_k \sin^2\theta_k &=& b_2
\nn \\ c_1 \cos\theta_1 \sin\theta_1 + c_2\cos\theta_2 \sin\theta_2 + \cdots + c_k
\cos\theta_k \sin\theta_k &=& 0 \end{eqnarray}

These equations, under the condition that $\sum_1^k c_i = b_1 + b_2$,  are equivalent to:
\begin{eqnarray} \sum_1^k c_i \cos (2\theta_i) &=& b_1-b_2 \nn
\\ \sum_1^k c_i \sin(2\theta_i) &=& 0 \end{eqnarray}

By defining $\displaystyle y_i = c_i \left[\begin{matrix} \cos 2\theta_i
\\ \sin 2\theta_i \end{matrix} \right]$, the above equations become the vector equation:

 \begin{equation} \sum_{i=1}^k y_i - \left[
\begin{matrix} b_1-b_2 \\ 0 \end{matrix} \right] = \left[ \begin{matrix} 0\\0 \end{matrix}\right]
\end{equation}

Using the ``tip-to-tail'' method of adding vectors, we can solve this equation by finding
a polygon with $k+1$ sides having lengths $c_1, c_2, \ldots, c_k, b_1-b_2$.  We must
place the side with length $b_1-b_2$ along the x-axis, but all of the angles between
sides are variable.  Such a polygon always exists for $k \geq 2$ provided the greatest
side length is less than or equal to the sum of the lengths of the remaining sides.  Since the
longest vector has either length $c_1$ or $b_1-b_2$, there is a solution when $c_1 \leq b_1$:

\begin{eqnarray*}
 c_1 &\leq& b_1 \\ 2c_1 &\leq& (b_1 - b_2) + (b_1 + b_2) \\ &=& (b_1-b_2) + \sum_{i=1}^k
c_i \\ c_1 &\leq& (b_1-b_2) + (c_2 + c_3 + \cdots + c_k) \\ \txt{and}\quad  b_1 - b_2 &\leq&
b_1 + b_2 \\ &=& c_1 + c_2 + \cdots + c_k
\end{eqnarray*}

\end{proof}

The following is the frame version of Theorem \ref{sums}.
%%%%%%%%%%%%%%%%%%%%%%%%
%    Corollary         %
%%%%%%%%%%%%%%%%%%%%%%%%
\begin{corollary}\label{fiframe}

Let $\hs$ be a Hilbert space with finite dimension $n$.  Let $B$ be a positive invertible operator on $\hs$ with
eigenvalues $b_1 \geq b_2 \geq \cdots \geq b_n > 0$, and for $k \geq n$ let $\{c_i\}_{i=1}^k$ be a sequence of
positive numbers with $c_1 \geq c_2 \geq \cdots \geq c_k$ such that $\sum_{i=1}^k c_i = \sum_{j=1}^n b_j$.  There
exists a frame for $\hs$ with $k$ vectors $\{x_i\}_{i=1}^k$ having frame operator $B$ and such that $\|x_i\| =
\sqrt{c_i}, 1 \leq i \leq k$, if and only if condition (\ref{ps}) from Theorem \ref{sums} is satisfied.
\end{corollary}

\begin{proof}
$B$ has a rank-one decomposition corresponding to $\{c_i\}_{i=1}^k$ if and only if there exist unit vectors
$\{y_i\}_{i=1}^k$ such that
$$B = \sum_{i=1}^k c_i(y_i \otimes y_i) = \sum_{i=1}^k \sqrt{c_i}y_i \otimes \sqrt{c_i}y_i$$
>From the expression (\ref{op}) of the frame operator, $\{x_i\}_{i=1}^k = \{\sqrt{c_i}y_i\}_{i=1}^k$ is a frame
with frame operator $B$ and $\|x_i\| = \sqrt{c_i}$.
\end{proof}

An immediate corollary is Theorem 5.1 from \cite{CKLT}, which finds the Fundamental Frame Inequality as the
necessary and sufficient condition for the existence of a tight frame with lengths of the vectors prescribed by
the sequence.  We state this result here.

\begin{corollary}\label{pete}\cite{CKLT} Let $\hs$ be a Hilbert space with finite dimension $n$, and
for some $k \geq n$, let $\{a_i\}_{i=1}^k$ be a sequence of positive numbers such that $a_1 \geq a_2 \geq \cdots
\geq a_k > 0$.  If $$ a_1^2 \leq \frac{1}{n} \sum_{i=1}^k a_i^2 \qquad \txt{(\textit{Fund. Frame Inequality})}$$
then there exist unit vectors $\{x_i\}_{i=1}^k$ such that the vectors $$ y_i = a_i x_i$$ form a tight frame for
$\hs$. The frame bound will be $$ \L = \frac{1}{n} \sum_{i=1}^k a_i^2$$

\end{corollary}

\begin{proof}
This case satisfies the hypotheses of Corollary \ref{fiframe} for $B = \L I$ and $c_i = a_i^2,\;i=1, \ldots, k$
since each $b_i = \L$, and $a_1^2 \leq \L$ implies the remaining inequalities.
\end{proof}

%%%%%%%%%%%%%%%%%%%%%%%%%%%%%%%%%%%%%%%%%%%%%%%%%%%%%%%%%%%%%%%%%
%     Section:  Infinite Dimensions                             %
%%%%%%%%%%%%%%%%%%%%%%%%%%%%%%%%%%%%%%%%%%%%%%%%%%%%%%%%%%%%%%%%%

\section{Infinite Dimensions}\label{two}
The main result, Theorem \ref{infinite}, of this section (or more particularly Corollary \ref{inframe}) is a
generalization of Theorem 5.4 from \cite{CKLT}, which we list here as Corollary \ref{frame}.  We will use
Proposition \ref{ident}, the operator-theoretic formulation of Corollary \ref{frame}, in the proof of Theorem
\ref{infinite}.

Corollary \ref{inframe} is also a generalization of Theorem 2 from \cite{DFK}.  This theorem states that, given
$B$ a positive operator with $\|B\|_{ess}>1$, there exists an ellipsoidal tight frame for $B$ with frame bound 1.
This is equivalent to the existence of a frame of unit vectors with frame operator $B^{-2}$, which is the special
case of Corollary \ref{inframe} with $\{c_i\}$ the constant sequence in which each $c_i = 1$.
%%%%%%%%%%%%%%%%%%%%
%   Theorem        %
%%%%%%%%%%%%%%%%%%%%
\begin{theorem}\label{infinite}  Let $B$ be a positive non-compact operator in $\bh$ for $\hs$ a
real or complex Hilbert space with infinite dimension.  If $\{c_i\}_{i=1}^{\infty}$ is a bounded
sequence of positive numbers with $\sup_i c_i < \|B\|_{ess}$ and $\sum_i c_i = \infty $, then
there is a sequence of rank-one projections $\{P_i\}_1^{\infty} \subset \bh$ such that
 $$ B = \sum_{i=1}^{\infty} c_i
P_i $$

\noi where convergence is in the SOT.   \end{theorem}

\begin{example} The condition that $\sup_i \{c_i\}$ be strictly less than $\|B\|_{ess}$ cannot be
dropped in general.  For instance, if $B=I$ the identity, and $$\{c_i\}_{i=1}^{\infty} = \{ \frac{2}{3}, 1, 1, 1,
\cdots \}$$ then there is no series $\{P_i\}$ of projections with $\sum_i c_i P_i = I$. Equivalently, there is no
sequence of vectors forming a Parseval frame in which the norms of the vectors $\{x_i\}$ are each $ \|x_i\| =
\sqrt{c_i}$.  If, in fact, such a frame did exist, the vectors $\{x_2, x_2, \ldots \}$ would necessarily form an
orthonormal set.  The assumption that the frame is Parseval then implies that $x_1$ is a unit vector, which is a
contradiction.
\end{example}

We will prove Theorem \ref{infinite} in a series of steps, using the finite dimensional Theorem
\ref{sums}.  The first step is Proposition \ref{ident}, which, as we previously stated, is a reformulation of Theorem 5.4 from \cite{CKLT}.   The result is
theirs, but our operator-theoretic proof is different from the one given in \cite{CKLT}.

\begin{proposition}\label{ident}
Let $\hs$ be an infinite-dimensional separable Hilbert space.  Let $\{c_i\}_{i=1}^{\infty}$ be a
sequence of numbers with $0 < c_i \leq 1$, and suppose $\sum_i c_i = \infty$ and $\sum_i (1-c_i) =
\infty$.  Then there is a sequence of rank-one projections $\{P_i\}_{i=1}^{\infty}$ such that
$$ I = \sum_{i=1}^{\infty} c_i P_i $$ with the sum converging in the strong operator topology.

\end{proposition}

An immediate corollary is the precise statement of Theorem 5.4 from \cite{CKLT}, which gives the existence of a Parseval frame for $\hs$ with norms of the frame
elements coming from the sequence.

\begin{corollary}\label{frame}\cite{CKLT}  If $\hs$ is an
infinite-dimensional real or complex Hilbert space, and if $\{a_i\}$ is a sequence of real numbers
with $0<a_i\leq 1$ such that $\sum_i a_i^2 = \infty$ and $\sum_i (1-a_i^2) = \infty$, then there is
a Parseval frame $\{x_i\}$ for $\hs$ such that $\|x_i\| = a_i$ for all $i$.
\end{corollary}

\begin{proof}[Proof of Proposition \ref{ident}]   For each $n \in \None$, let $s(n)$ denote the
smallest integer such that $$c_1 + c_2 + \cdots + c_{s(n)} > n$$  (Such a value exists for every
$n$ since $\sum_i c_i = \infty$.)  Denote what we will call the \textit{residual} by $$r(n) = c_1 +
c_2 + \cdots + c_{s(n)} - n$$ and denote the \textit{integer gap} by $$g(n) = s(n)-n$$

Then $0 < r(n) \leq 1$, and clearly $g(n) \leq g(n+1)$ for all $n \in \None$.  Denote $\D_i =
1-c_i$, and let $d(n) = \D_1 + \D_2 + \cdots + \D_n$.  We have
\begin{eqnarray*}  g(n) = s(n)-n &>& s(n)-(c_1 + \cdots c_{s(n)}) \\ &=& d(s(n)) \\ &\geq& d(n)
\end{eqnarray*}

By hypothesis, $\sum_i \D_i = \infty$, so $d(n) \rightarrow \infty$ and thus $g(n) \rightarrow \infty$. We can
therefore inductively choose an increasing sequence $\{n_i\}_{i=1}^{\infty}$ of natural numbers such that for all
$i \geq 2$,
\begin{equation} \label{a}
g(n_i) \geq g(n_{i-1}) + 1  \end{equation} and
\begin{equation} \label{b}
d\left(s(n_{i-1})+n_i-n_{i-1}+1 \right) - d\left( s(n_{i-1}) \right) > 2
\end{equation}

Note that condition \ref{a} implies that $s(n_i)-s(n_{i-1}) \geq n_i - n_{i-1}+1$, and condition
\ref{b} implies that:

\begin{equation}\label{c}
 c_{(s(n_{i-1}+1))} + \cdots + c_{(s(n_{i-1}) + n_i - n_{i-1} + 1)} < n_i - n_{i-1}-1
 \end{equation}

Let $\{e_i\}_1^{\infty}$ be an orthonormal basis for $\hs$, and let $E_i = e_i \otimes e_i$ for all $i \in \None$.
Let $$ B_1 = E_1 + E_2 + \cdots + E_{n_1} + r(n_1)E_{n_1+1}$$  and for $i \geq 2$, let $$ B_i =
\left((1-r(n_{i-1})\right) E_{n_{i-1}+1} + \left[ E_{n_{i-1}+2} + \cdots + E_{n_i} \right] + r(n_i)E_{n_i + 1}
$$ This gives $$ \tx{trace}(B_1) = c_1 + \cdots + c_{s(n_1)}$$ and for $i \geq 2$,
\begin{eqnarray*}
 \tx{trace}(B_i) &=& [1-r(n_{i-1})] + [n_i - n_{i-1} - 1] + r(n_i) \\
 &=& [r(n_i) + n_i] - [r(n_{i-1})+n_{i-1}] \\
 &=& \sum_{j=1}^{s(n_i)} c_j - \sum_{j=1}^{s(n_{i-1})} c_j \\
 &=& \sum_{j=s(n_{i-1})+1}^{s(n_i)} c_j
 \end{eqnarray*}

 We have rank$(B_1) = n_1 + 1$, and for $i \geq 2$, we have rank$(B_i) \leq n_i - n_{i-1}+1$.  When
$i \geq 2$, at least the first $n_i - n_{i-1}-1$ (in decreasing order) of the eigenvalues of $B_i$ are $1$.  (That
is, all but the last two eigenvalues of $B_i$ are 1.) This, together with the condition (\ref{c}), implies that
$B_i$ and the sequence $\{c_{s(n_{i-1})+1}, \ldots , c_{s(n_i)}\}$ satisfy the hypotheses of Theorem \ref{sums}.
To see this, for $1 \leq p \leq n_i - n_{i-1}-1$, the inequalities (\ref{ps}) in Theorem \ref{sums} are trivially
satisfied.  For $p=n_i - n_{i-1}$, the inequality (\ref{c}) implies that the sum of the first $n_i - n_{i-1}$ of
the terms in $\{c_{s(n_i-1)+1}, \ldots, c_{s(n_i)} \}$ is less than $n_i - n_{i-1}-1$, which is less than or equal
to the sum of the first $n_i - n_{i-1}$ eigenvalues of $B_i$.

 Therefore, for $i \geq 2$, there exist rank-one projections $\{P_{s(n_{i-1})+1}, \ldots ,
P_{s(n_i)}\}$ in the range of $B_i$ such that
$$ B_i = c_{s(n_{i-1})+1}P_{s(n_{i-1})+1} + \cdots + c_{s(n_i)}P_{s(n_i)} $$

For the special case $B_1$, since there is at most one eigenvalue different from 1, the hypotheses of Theorem
\ref{sums} are clearly satisfied, and so there exist rank-one projections $P_1, \ldots P_{s(n_1)}$
such that $$B_1 = \sum_{j=1}^{s(n_1)} c_j P_j $$

For all $i \in \None$,
$$B_1 + \cdots + B_i = E_1 + \cdots + E_{n_i} + r(n_i)E_{n_i + 1}$$
Since $\{r(n_i)\}$ is bounded, the sum $\sum_{i=1}^{\infty} B_i$ converges strongly to the
identity $I$.  It follows that $\sum_{i=1}^{\infty} c_i P_i$ converges strongly to $I$.
\end{proof}

\begin{corollary}\label{subspace} Let $E$ be an infinite rank projection in $\bh$, and let $\{c_i\}$
be a sequence of real numbers with $0<c_i\leq 1$ such that $\sum_i c_i = \infty$ and $\sum_i
(1-c_i) = \infty$.  Then $E$ has a rank-one decomposition corresponding to $\{c_i\}$:
$$E = \sum_{i=1}^{\infty} c_i P_i $$
\end{corollary}
\begin{proof}  Apply Proposition \ref{ident} to the subspace $E\hs$ and the result follows
immediately.
\end{proof}

In \cite{DFK}, it was proven (Theorem 2) that if $B$ is a positive operator with essential norm strictly greater
that 1, then $B$ has a projection decomposition.  That is, $B= \sum_i P_i$, where $\{P_i\}$ are self-adjoint
projections and convergence is in the strong operator topology.  Clearly $\{P_i\}$ can be taken to
be rank-one.  We need here a stronger result, that each of the projections can be required to have
infinite rank.

\begin{proposition}\label{proj}

Let $B$ be a positive operator in $\bh$ for $\hs$ a real or complex separable Hilbert space with
infinite dimension, and suppose $\|B\|_{ess} > 1$.  Then $B$ has an infinite-rank projection
decomposition.
\end{proposition}

\begin{proof}
The hypothesis that $\|B\|_{ess}>1$ implies there exists a sequence $\{Q_i \}_{i=1}^{\infty}$ of mutually
orthogonal infinite-rank projections commuting with $B$ such that $\|B|_{Q_i\hs} \|_{ess} > 1$ for all $i$, and
with $\sum_i Q_i = I$.  Let $A_i = B|_{Q_i\hs}$.  By Theorem 2 in \cite{DFK}, each $A_i$ has a projection
decomposition: $$ A_i = \sum_{j=1}^{\infty} E_{ij} $$ Here we can assume without loss of generality that for each
$i$, $\{E_{i,j}\}_{j=1}^{\infty}$ is an infinite sequence of nonzero projections because $A_i$ is not finite rank.
Note that for each $j$, the projections $\{ E_{1,j}, E_{2,j}, \ldots, \}$ are mutually orthogonal.  Let $$E_j =
E_{1,j} + E_{2,j} + \cdots  $$ with the sum converging strongly.  Then $E_j$ is an infinite-rank projection, and
we have $\sum_{j=1}^{\infty} E_j = B$, as required.

\end{proof}

\begin{remark}\label{part} If $\mathcal{E}$ is a countably infinite set of positive numbers,
then it is well-known that $\sum_{\L \in \mathcal{E}} \L$ makes sense as an extended real
number, independent of enumeration.  It is also clear that if $\sum_{\L \in \mathcal{E}} \L =
\infty$, then there exists a partition $\{\mathcal{E}_i \}_{i=1}^{\infty}$ of $\mathcal{E}$ into
infinitely many subsets such that for each $i$, $\sum_{\L \in \mathcal{E}_i} \L = \infty$.
Similarly, if $\{c_i\}_{i=1}^{\infty}$ is a sequence of positive numbers such that
$\sum_i c_i = \infty$, there exists a partition of the sequence into infinite subsequences
$\Lambda_i = \{c_{i,j} \}_{j=1}^{\infty}$ such that $\sum_{j=1}^{\infty} a_{i,j} = \infty$ for each
$i$. We will use this partition in the proof of the next theorem. \end{remark}

The reader will note that Theorem \ref{oper} is an alternate restatement of
Theorem \ref{infinite}.  We include it because it has a simpler form and has
independent interest.  We prove it first and then obtain Theorem \ref{infinite} by a
simple scaling argument.

\begin{theorem}\label{oper}   Let $B$ be a positive operator in $\bh$ for $\hs$ with $\|B\|_{ess}>1$.
Let $\{c_i\}_{i=1}^{\infty}$ be any sequence of numbers with $0<c_i \leq 1$ such that $\sum_i c_i = \infty$. Then
there exists a sequence of rank-one projections $\{P_i\}_{i=1}^{\infty}$ such that
$$ B = \sum_{i=1}^{\infty} c_i P_i $$
\end{theorem}

\begin{proof}
Choose $\A,\; 0<\A < 1$ such that $\|\A B\|_{ess} > 1$. Using
Proposition \ref{proj}, write $\A B = \sum_n Q_n$, a projection decomposition with each $Q_n$ having infinite rank. Next, let $\G_i = \A c_i$ for all $i$, and partition the sequence $\{\G_i\}$ into subsequences $\Lambda_n$ as described in Remark \ref{part}.  Note that for each $n$, $\sup_j \{\G_{n,j} \} < 1$.   
Apply Corollary \ref{subspace} to each operator $Q_n$ with the subsequence $\Lambda_n$, obtaining rank-one
projections $\{F_{n,j} \}_{j=1}^{\infty}$ with $$Q_n = \sum_{j=1}^{\infty} \G_{n,j} F_{n,j}$$ Scaling by $\A$, 
we have the rank-one decomposition of $B$: $$B = \frac{1}{\A} \sum_{n=1}^{\infty} \sum_{j=1}^{\infty}
\G_{n,j}F_{n,j} = \sum_{n=1}^{\infty} \sum_{j=1}^{\infty}
c_{n,j}F_{n,j} $$

\end{proof}

\begin{proof}[Proof of Theorem \ref{infinite}]  We apply Theorem \ref{oper}.  Given $B$ non-compact, we have $m := \|B\|_{ess} >
0$ and $\sup_i c_i < m$.   Take $\A = \frac{1}{m-\E}$ for some $\E > 0$, so that we have $\|\A B
\|_{ess} > 1$, and let $\G_i = \A c_i$, which gives $\sup_i \G_i < 1$.  By Theorem \ref{oper}, there
is a rank-one decomposition for $\A B$ corresponding to the sequence $\{\G_i\}$:
$$\A B = \sum_{i=1}^{\infty} \G_i P_i$$
Scaling the above by $m-\E$ gives the final result.
\end{proof}

The following is the frame-theoretic version of Theorem \ref{infinite}.  It is now clear that
 Theorem 5.4 from \cite{CKLT} is the special case of Corollary \ref{inframe} in which $B$ is the identity operator, and Theorem 2 from \cite{DFK} is the case where $\{c_i\}$ is a constant sequence.

\begin{corollary}\label{inframe} Let $\hs$ be a real or complex separable Hilbert space of
infinite dimension, and let $B$ be any bounded positive invertible operator
on $\hs$.  Let $\{a_i\}$ be an arbitrary sequence of positive numbers with  infinite sum and with $\sup_i \{a_i^2\}  < \|B\|_{ess}$.  Then there is a frame $\{z_i\}$ for $\hs$ with frame operator $B$
such that $\|z_i\| = a_i \;\txt{for all} i$.
\end{corollary}
\begin{proof}
Let $c_i = a_i^2$, and let $\A = \sup_i \{a_i^2\} = \sup_i \{c_i\}$.  Then $\A < \|B\|_{ess}$, so by Theorem \ref{infinite}, $B$ has a rank-one decomposition corresponding to $\{c_i\}$.  Let $\{x_i\}$ be unit vectors in the range of the rank-one operators of this decomposition, and let $z_i = \sqrt{c_i}x_i$ for all $i$.   $$ B = \sum_i c_i (x_i \otimes x_i) = \sum_i z_i \otimes z_i$$  By the expression (\ref{op}) of the frame operator, $\{z_i\}$ is the desired frame.
\end{proof}

%%%%%%%%%%%%%%%%%%%%%%%%%%%%%%% %% Bibliography %% %%%%%%%%%%%%%%%%%%%%%%%%%%%%%%%
\bibliographystyle{plain}
%%%%%%%%%%%%%%%%%%%%%%%%%%%%%%%%%%%%%%%%%%%%%%%%%%%%%%%%%
% \nocite{*}
%   \bibliography{cfrefs}
%%%%%%%%%%%%%%%%%%%%%%%%%%%%%%%%%%%%%%%%%%%%%%%%%%%%%%%%%

\end{document}